\documentclass{amsart}
\usepackage{amsthm}
\usepackage{amssymb}

\newtheorem{Theorem}{Theorem} [section]
\newtheorem{Lemma}[Theorem]{Lemma}
\newtheorem{Condition}[Theorem]{Condition}

\theoremstyle{definition}
\newtheorem{Definition}[Theorem]{Definition}

\theoremstyle{remark}
\newtheorem{Remark}[Theorem]{Remark}
\newtheorem*{explanation}{Explanation}
\newtheorem{Example}[Theorem]{Example}

\numberwithin{equation}{section}

\def	\C 	{{\mathbb C}}
\def	\R 	{{\mathbb R}}
\def	\Z 	{{\mathbb Z}}
\def	\g 	{{\mathfrak g}}
\def	\h 	{{\mathfrak h}}
\def	\cut 	{{\text{cut}}}
\def	\ssminus 	{\smallsetminus}
\def	\Nu	{{\mathcal V}}

\def	\Stab	{{\operatorname{Stab}}}
\def	\inv	{^{-1}}
\def	\ol	{\overline}
\def	\calA	{{\mathcal A}}
\def	\tZ	{\tilde{Z}}
\def	\tY	{\tilde{Y}}
\def	\tX	{\tilde{X}}

\def	\tN	{\tilde{N}}
\def	\tPsi	{\tilde{\Psi}}

\newcommand {\comment}[1] {}
\newcommand{\printname}[1]{}
\newcommand{\labell}[1] {\label{#1}\printname{#1}}

\begin{document}

\title[Compact and non-compact equivariant cobordisms]
{The relation between compact and non-compact equivariant cobordisms}

\author{Viktor L. Ginzburg}
\address{Department of Mathematics, University of California
at Santa Cruz, Santa Cruz, CA 95064}
\email{ginzburg@count.UCSC.EDU}

\author{Victor L. Guillemin}
\address{Department of Mathematics, Massachusetts Institute of Technology,
Cambridge, MA 02139}
\email{vwg@math.mit.edu}

\author{Yael Karshon}
\address{Institute of Mathematics, The Hebrew University of Jerusalem,
Jerusalem 91904, Israel}
\email{karshon@math.huji.ac.il}

\thanks{The authors are supported in part by the NSF and by the BSF}

\date{November 30, 1998}

\begin{abstract}
We show that the theory of stable complex $G$-cobordisms, for a torus
$G$, is embedded into the theory of stable complex $G$-cobordisms of
not necessarily compact manifolds equipped with proper abstract moment
maps. Thus the introduction of such non-compact cobordisms in the stable
complex $G$-cobordism theory does not lead to new relations.
\end{abstract}

\maketitle

\section{Proper abstract moment maps.}
\labell{sec:intro}

The main objective of this paper is to show that geometric
$G$-equivariant cobordism theory is embedded in a similar theory
for non-compact manifolds.  Recall that two compact $G$-manifolds,
where $G$ is a compact Lie group, are said to be cobordant if their 
disjoint union is a boundary of a compact $G$-manifold. Often, the manifolds
are assumed to carry an additional structure preserved by the action,
and this structure is assumed to extend over the cobordism. The 
structures important for our present purpose are the orientation
or/and tangential stable complex structure. (For definitions, see 
Appendix \ref{sec:stable-complex} of this paper or 
Chapter 28 by G. Comeza\~{n}a in \cite{may}.) The cobordism classes of
$G$-manifolds are then referred to as (geometric) oriented or stable complex
$G$-cobordisms. We restrict our attention to the case where $G$ is a 
torus. (The adjective ``geometric'', omitted from now on, is used here
to distinguish the cobordism theory we consider from its
homotopy theoretic counterpart. See, e.g.,  Chapter 15 by 
S. R. Costenoble in \cite{may}.)

In \cite{K:cobordism} we introduced a cobordism theory whose objects
are non-compact $G$-manifolds with yet an additional structure called
a proper abstract moment map. (We will recall its definition below.) The 
reason for considering non-compact manifolds is that this allows one 
to obtain a simple form of the linearization
theorem (see \cite{GGK1} and \cite{K:cobordism}). In its non-compact
version proved in \cite{K:cobordism}, the linearization theorem claims
that under certain natural hypotheses every $G$-manifold is cobordant
to the normal bundle to its fixed point set with a suitable proper 
abstract moment map. The addition of non-compact manifolds to a 
cobordism theory could, however, create a problem. Namely, as a result,
all compact manifolds might then become cobordant to each other and so 
to the empty set. So to say,
adding non-compact manifolds has a trivializing effect on the cobordism
theory. For example, in the cobordism theory of $G$-manifolds equipped
with proper real-valued $G$-invariant functions every compact manifold 
is cobordant to zero.

In this paper we show that this does not happen for cobordisms with
proper abstract moment maps: the theory of stable complex $G$-cobordisms, for
a torus $G$, is embedded into the theory of stable complex $G$-cobordisms of
not necessarily compact manifolds equipped with proper abstract moment
maps. Thus the introduction of such new objects in stable complex $G$-cobordism
theory does not lead to new relations.

Let us now recall some definitions. In what follows $G$ is a torus
and all manifolds are assumed to be equipped with a
$G$-action.

\begin{Definition} 
An \emph{abstract moment map} is a $G$-equivariant map 
$ \Psi \colon M \to \g^*$ such that for any subgroup $H \subset G$,
the composition of $\Psi$ with the natural projection $\g^* \to \h^*$
is locally constant on the set of points fixed by $H$.
\end{Definition}

Examples and a detailed discussion of this notion can be found
in \cite{K:cobordism} and \cite{GGK-assign}. Here we only note
that on a compact manifold $M$, the identically zero function, 
$\Psi\equiv 0$, is a proper abstract moment map 

\begin{Definition} \labell{cobordism}
Let $\Psi_j \colon M_j \to \g^*$, $j=1,2$, 
be proper abstract moment maps on $G$-manifolds. 
A \emph{proper cobordism} between $M_1$ and $M_2$ is a 
$G$-manifold with boundary $M$, a proper abstract moment map
$\Psi \colon M \to \g^*$, and an equivariant diffeomorphism
\begin{equation} \label{del M}
 \partial M \cong M_1 \sqcup M_2
\end{equation}
which carries $\Psi$ to $\Psi_1 \sqcup \Psi_2$.
\end{Definition}

When, in addition, the manifolds are equipped with an extra structure
(such as a tangential $G$-equivariant stable complex structure),
this structure is assumed to extend over the
cobording manifold $M$ in the standard way. (See, e.g., \cite{stong}.)
The theory of compact $G$-cobordisms is mapped into its non-compact
counterpart for $G$-manifolds with proper abstract moment maps by
equipping every compact manifold with the identically zero abstract moment
map.  As has been mentioned above, for stable complex cobordisms this map
is one-to-one. More explicitly, in Section \ref{sec:compact-versus-proper}
we prove the following

\begin{Theorem}
\labell{thm:intro}
Let $M_1$ and $M_2$ be compact stable complex $G$-manifolds which are
properly cobordant when equipped with identically zero abstract moment maps. 
Then $M_1$ and $M_2$ are compactly cobordant.
\end{Theorem}

\begin{Remark} \labell{any maps}
In this theorem, the zero abstract moment maps can be
replaced by any abstract moment maps on $M_1$ and $M_2$.
\end{Remark}

\begin{Remark}
Manifolds with stable complex structures are automatically
orientable but the orientation is not canonical. 
An orientation is usually fixed in addition to a stable complex 
structure as the data giving the stable complex cobordism class. 
This is necessary for Chern numbers to be well defined.
In the present paper we treat both structures independently.
(See Appendix \ref{sec:stable-complex} below for references and
a detailed discussion.) Theorem \ref{thm:intro} still holds when all
the manifolds are oriented in addition to being stable complex.
\end{Remark}

If we allow orbifold cobordisms (see, e.g., \cite{Dr}), the theorem
(without the requirement that $M_1$ and $M_2$ are stable complex)
can be easily proved by using the \emph{Lerman cutting} of the 
cobording manifold. Namely, by cutting the cobording manifold 
with respect to a Delzant polytope we obtain a compact
orbifold. Choose the polytope so that its interior contains
the origin or, in the case of Remark \ref{any maps},
make it large enough so that its
interior contains the moment map image of the (compact) boundary
$M_1\sqcup M_2$. 
Then the cut orbifold gives a compact orbifold cobordism between $M_1$
and $M_2$. 
In Section \ref{sec:cutting}
we provide the details of this cutting argument. 

Unfortunately, in general, the resulting cobordism is only an
orbifold cobordism.   In Section \ref{sec:surgery}
we describe a surgery, dating back to a result of
Gusein-Zade \cite{GZ:short}, which gets rid of the orbifold 
singularities of the cut space. The proof of the theorem, in
a form more general than given above, is finished in Section 
\ref{sec:compact-versus-proper}.

Finally, we would like to note that Gusein-Zade in \cite{GZ:short} 
described the stable complex cobordism group of compact $S^1$-manifolds.

\section{Lerman cutting}
\labell{sec:cutting}

Let $M$ be a manifold with an action of a torus $G$
and an abstract moment map $ \Psi \colon M \to \g^*$.
Let $S^1 \subseteq G$ be a subcircle, generated by a Lie algebra element
$\eta \in \g$; it acts on $M$ with an abstract moment map 
$\Psi^\eta(m) = \langle \Psi(m), \eta \rangle.$
Let $a \in \R$ be a regular value for $\Psi^\eta$.
As a topological space, the \emph{cut space} $M_\cut$ 
is the quotient $\{ \Psi^\eta \geq a \} / \sim$, where 
on the level $ (\Psi^\eta) \inv (a) $ we have $m \sim m'$ 
if and only if $m$ and $m'$ 
in the same $S^1$-orbit and on the open part $\{ \Psi^\eta > a \}$ 
the relation $\sim$ is trivial: $m \sim m'$ if and only if $m=m'$.
The $G$-action and the abstract moment map on $M$ descend to $M_\cut$.

To define the cut space as a $C^\infty$ manifold (or orbifold)
we need to treat it more carefully.
Consider the product $M \times \C$ with the diagonal $S^1$-action
and the abstract moment map $\varphi(m,z) = \Psi^\eta(m) - |z|^2$.  
The \emph{cut space} with respect to $\eta$ at the value $a$ 
is the quotient 
\begin{equation} \labell{cut}
   M_\cut = Z / S^1, \quad \quad \text{where} \quad  Z = \{ \varphi = a \}.
\end{equation}
It is easy to check that if $a$ is a regular value for $\Psi^\eta$,
it is also a regular value for $\varphi$. Let us assume that this is 
the case. Then $Z$ is a manifold.
The diagonal action of $S^1$ on the level set $Z=(\varphi) \inv (a)$ has 
finite
stabilizers; we proved this in \cite[Lemma 7.1]{K:cobordism} for any
regular level set of an abstract moment map. Therefore, the 
cut space is an orbifold.

The product $M \times \C$ inherits from $M$ 
a left $G$-action and an abstract moment map on $M$,
which descends to $M_\cut$. 
A group action on $M$ which commutes with $G$ and preserves $\Psi$
also descends to $M_\cut$.
Similarly, if $M$ is equivariantly stable complex or oriented, 
so is $M_\cut$.

\begin{Remark}
The cutting construction described above was invented by  Lerman
\cite{L:cut} for symplectic manifolds with Hamiltonian group actions.
An alternative description, perhaps more familiar in topology,
is that $M_\cut$ is obtained by gluing the disk (orbi-)bundle 
$(\Psi^\eta) \inv(a) \times_{S^1} D^2$ to $\{ \Psi^\eta \geq a \}$. 
However, we find this definition based on the surgery to be less 
convenient 
to work with than the original Lerman construction. The reason is that 
carrying out the surgery requires a choice 
(of a collar around $ ( \Psi^\eta ) \inv (a)$),
and this choice makes it more difficult to show that the cut space 
naturally inherits various structures from $M$.
\end{Remark}

The cutting construction provides a way to turn a non-compact
cobordism with a proper abstract moment map into  compact one
at the cost of introducing singularities.

\begin{Theorem} \labell{orbifold}
Let $M_1$ and $M_2$ be two compact $G$-manifolds with abstract moment 
maps. Suppose that $M_1$ and $M_2$ are properly cobordant in the sense 
of Definition  \ref{cobordism}. Then these manifolds are also cobordant 
through a \emph{compact} orbifold with a $G$-action and an abstract 
moment map.

Moreover, if the manifolds and the proper cobordism possess one or more 
of the following structures, then so does the compact orbifold cobordism:
\begin{enumerate}
\item
an action of another group which commutes with the $G$-action
and preserves the $G$-moment map;
\item
an equivariant stable complex structure;
\item
an orientation.
\end{enumerate}
\end{Theorem}

\begin{proof}
Let $M$ be a (possibly non-compact) cobording manifold between $M_1$
and $M_2$. Pick elements $\eta_1, \ldots, \eta_r$ of $\g$ whose 
positive span is $\g$. Then for any real numbers 
$a_1, \ldots, a_r$, the subset 
\begin{equation} \labell{polytope}
 \bigcap \{ \eta_j \geq a_j \}
\end{equation}
of $\g^*$ is a compact polytope.
We can choose all $\eta_j$ to be integral, i.e., to generate
circle subgroups of $G$.
Furthermore, we can choose  $a_1$ to be a regular value of 
$\Psi^{\eta_1}$ which is small enough so that $\Psi$ sends the 
boundary of $M$ 
into the open half-space $\{ \eta_1 > a_1 \}$ of $\g^*$.
This is possible because, by assumption, $\partial M$ is compact.
We cut $M$ with respect to $\eta_1$ at the value $a_1$,
as described above.
This produces a cobording orbifold, $M'$, with a $G$-action
and a proper abstract moment map whose image 
is contained in the closed half-space $\{ \eta_1 \geq a_1 \}$. 
If $M$ possesses another group action, an equivariant stable complex
structure, and/or an orientation, $M'$ inherits these structures.

We proceed by induction, each time cutting with respect to $\eta_j$ 
at a regular value $a_j$ which is small enough so that the boundary is sent 
into $\{ \eta_j > a_j \}$. Once we reach $j=r$, we obtain a cobording 
orbifold whose image is contained in the polytope \eqref{polytope}. 
Since the abstract moment map is proper and \eqref{polytope} is compact, 
this cobording orbifold is also compact.
\end{proof}

\begin{Remark}
One can cut the cobordism $M$ simultaneously over all the faces of the
polytope as in \cite{convexity}. However, we prefer to perform the cutting 
sequentially, because this allows us to work with $S^1$-quotients only,
and their singularities can be resolved. 
\end{Remark}

\begin{Remark}
We did not recall the definition of an orbifold, a group action on an
orbifold, and a stable complex structure on an orbifold.
(See \cite{orbi}.)
In the context of this section, it suffices to know that if $S^1$ 
acts on $Z$ with finite stabilizers, then $Z/S^1$ is an orbifold, 
and that a group action on $Z$ which commutes with the $S^1$-action
and an equivariant stable complex structure on $Z$ both descend to the 
orbifold $Z/S^1$.
The statements and proofs in Section 
\ref{sec:compact-versus-proper} (including the proof of Theorem 
\ref{thm:intro}) are complete and accurate without and, in fact,
independent of explicit definitions of these structures on orbifolds. 
\end{Remark}

\begin{Remark} \labell{quasi free}
If $G=S^1$ and the action on the cobordism $Z$ is quasi-free, i.e.,
with stabilizers either $\{ 1 \}$ or $S^1$, 
then the compact cobordism that we get
is an actual manifold, without any orbifold singularities, and the action 
on it is still quasi-free.
\end{Remark}

\section{Gusein-Zade surgery}
\labell{sec:surgery}

The orbifold singularities of the Lerman cut space result from taking
in \eqref{cut} the quotient 
by the $S^1$-action, which has finite, but perhaps non-trivial,
stabilizers. In this section we describe a surgery which gets rid of the
singularities by replacing the action on $Z$ by a free action. This
construction has been used repeatedly by various mathematicians; 
the earliest reference that we found is a 1971 paper by Gusein-Zade,
\cite{GZ:short}.

Recall that the orbit type stratification of a manifold with respect
to an $S^1$-action is the decomposition of the manifold into connected 
components of the sets $\{$\emph{points whose stabilizer is $\Gamma$}$\}$
for subgroups $\Gamma$ of $S^1$.  
The strata are partially ordered; $X \leq X'$ if and only if
$X$ is contained in the closure of $X'$. A stratum is said to be 
\emph{minimal} if it is closed or, equivalently, if it does not contain
other strata in its closure.

The main result of this section is 

\begin{Lemma} \labell{desingularize}
Let $Z$ be an $S^1$-equivariant stable complex manifold 
that has no fixed points. Let $Y \subset Z$ be the set of points where
the action is not free.
Then one can obtain a manifold with a free $S^1$-action 
by performing an equivariant surgery on $Z$ along $Y$. 

The surgery can be performed so as to preserve one or more 
of the following additional structures:
\begin{itemize}
\item
an orientation 
\item
a $G$-action which commutes with the circle action, where $G$
is a compact Lie group, and a $G \times S^1$-equivariant 
stable complex structure which lifts the given $S^1$-equivariant
stable complex structure.
\item
a proper abstract moment map.
\end{itemize}
The resulting manifold, $Z'$, is equivariantly cobordant to $Z$ with 
the same additional structures.
\end{Lemma} 

\begin{Remark}
The condition that the manifold be equivariantly stable complex
can be replaced by an assumption on the stabilizers. See Section 
\ref{sec:null}.
\end{Remark}

Before proving the statement of Lemma \ref{desingularize}, 
we pause to explain it: 

\begin{explanation}
We will construct a manifold $Z'$ with a free $S^1$ action,
and a closed invariant subset $Y' \subset Z'$,
and an equivariant diffeomorphism 
\begin{equation} \labell{eq1a} 
 Z \ssminus Y \cong Z' \ssminus Y'.
\end{equation}

We will have an $S^1$-manifold $\tZ$ with boundary,
a closed invariant subset $\tY \subset \tZ$
which is a locally finite union of closed invariant submanifolds
transverse to the boundary, and an equivariant diffeomorphism 
\begin{equation}\labell{eq2a}
 \partial \tZ \cong Z \sqcup Z'
\end{equation}
which carries $\partial \tY$ to $Y \sqcup Y'$.
On the open dense complement of these subsets, 
the cobordism will
be trivial, i.e., we will have an equivariant diffeomorphism
\begin{equation}\labell{eq3a}
 \tZ \ssminus \tY \cong V \times [0,1]
\end{equation}
and isomorphisms 
$V \cong Z \ssminus Y \stackrel{\eqref{eq1a}}{\cong} Z' \ssminus Y'$.

Often one starts with $Z$ which is compact and demands that
$Z'$ and $\tZ$ be compact. We do not assume compactness; instead, 
we will have the following condition:
for any subset $A \subset V \cong Z \ssminus Y$
whose closure in $Z$ is compact, 
the closure of $A \times [0,1]$ in $\tZ$
(see \eqref{eq3a}) is also compact.
Informally speaking, $Z'$ and $\tZ$ have no holes. 
In particular, if $Z$ is compact, so are $Z'$ and $\tZ$.

If $Z$ is oriented, then so are $Z'$ and $\tZ$;
\eqref{eq1a} respects orientations, and \eqref{eq2a}
carries the boundary orientation of $\partial \tZ$
to the given orientation on $Z$ and the opposite orientation on $Z'$.
If $Z$ has a $G$-action which commutes with the $S^1$-action, 
with $G$ a compact Lie group, then so do $Z'$ and $\tZ$, 
and \eqref{eq1a} and \eqref{eq2a} are $(G \times S^1)$-equivariant. 
If $Z$ has a proper abstract moment map, then so do $Z'$ and $\tZ$,
\eqref{eq2a} respects these maps, and \eqref{eq1a} respects these maps 
outside neighborhoods of $Y$ and $Y'$ which can be
pre-chosen to be arbitrarily small.
An equivariant stable complex structure on $Z$ induces such
structures on $Z'$ and $\tZ$, and \eqref{eq1a} and \eqref{eq2a} lift
to isomorphisms of equivariant stable complex manifolds (see Appendix
\ref{sec:stable-complex}).
\end{explanation}

We prove Lemma \ref{desingularize} by inductively replacing the non-free
orbit type strata by strata with smaller stabilizers.
The following lemma provides the induction step:

\begin{Lemma} \labell{step}
Let $Z$ be a manifold with an $S^1$ action that has no fixed points.
Let $X \subset Z$ be a minimal orbit type stratum.
Suppose that the normal bundle, $\Nu$, of $X$ in $Z$ 
admits an equivariant fiberwise complex structure.
Then by performing an equivariant surgery on $Z$ along $X$
one can replace $X$ by strata with smaller stabilizers.

The surgery can be carried out so as to preserve one or more of the
following additional structures:
\begin{itemize}
\item
an orientation,
\item
a $G$-action which commutes with the circle action,
where $G$ is a compact Lie group, preserving the complex structure on $\Nu$,
\item
a proper abstract moment map,
\item
an equivariant stable complex structure.
\end{itemize}
The resulting manifold, $Z'$, is equivariantly cobordant to $Z$ with the
same additional structures.
\end{Lemma} 

\begin{explanation}
An explanation of Lemma \ref{step} 
is completely analogous to that of Lemma \ref{desingularize}
with $Y$, $Y'$, and $\tY$ replaced by $X$, $X'$, and $\tX$,
except that the action on $Z'$ need not be free.
All the points in $X'$ will have stabilizers of orders smaller
than the order of the stabilizer of $X$.
\end{explanation}

Lemma \ref{step} implies Lemma \ref{desingularize} by induction:

\begin{proof}[Proof of Lemma \ref{desingularize}]
Let $X_1$, $X_2$, $\ldots$ be all the minimal strata in $Z$.
Because $Z$ is equivariantly stable complex, the normal bundles of $X_j$
in $Z$ are equivariantly complex.
The surgery of Lemma \ref{step} can be applied to all of them
simultaneously, because they are closed and disjoint.
Denote the resulting manifold by $Z_1$. 
Proceed inductively to construct a sequence $Z_1, Z_2, \ldots$,
where each $Z_{n+1}$ is obtained by applying the surgery of
Lemma \ref{step} to all the minimal strata of $Z_n$. 

To explain the convergence of this process, exhaust $Z$ by open subsets
whose closures are compact: $Z = \cup_j C^j$.
For each fixed $j$, the sequence of surgeries give open subsets
$C_1^j \subseteq Z_1$, $C_2^j \subseteq Z_2$, etc., each with a
compact closure, and cobordisms between $C_1^j$ and $C_2^j$,
$C_2^j$ and $C_3^j$, etc. At each stage, the order of the largest
stabilizer in $C_{m+1}^j$ is strictly smaller than that in $C_m^j$. 
Therefore, after finitely many steps, the action on $C_m^j$ will be
free, and we will get $C_m^j = C_{m+1}^j = \ldots$, which we denote
$C_\infty^j$, with an infinite sequence of trivial cobordisms.
By adding a copy of $C_\infty^j$ at the very end, we can ``close up"
the infinite sequence of cobordisms to obtain a cobordism between 
$C^j$ and $C_\infty^j$.

The stabilization might occur at a different place $m$ for different $j$'s. 
However, the $C_\infty^j$'s fit together into a well defined limit, 
$Z_\infty$, and the cobordisms between $C^j$ and $C_\infty^j$ fit together
and form a cobordism between $Z$ and $Z_\infty$.

The additional requirements in Lemma \ref{desingularize}
follow immediately from those in Lemma \ref{step}.
\end{proof}

The proof of Lemma \ref{step} will rely on the following key observation:

\begin{Lemma} \labell{extend}
Let $Z$ be a manifold with an $S^1$-action,
$X \subseteq Z$ an orbit type stratum with finite stabilizer $\Gamma$,
and $\Nu$ the normal bundle to $X$ in $Z$ with the induced $S^1$-action.
Suppose that $\Nu$ admits an equivariant fiberwise complex structure.
Then the fiberwise action of $\Gamma$ on $\Nu$ 
extends to a fiberwise action of $S^1$, 
which commutes with the original $S^1$-action,
and for which the stabilizer of any point outside of the zero section
has fewer elements than $\Gamma$.

If a compact Lie group $G$ acts on $\Nu$, preserves its complex
structure, and commutes with the $S^1$-action, we can choose the fiberwise
$S^1$-action on $\Nu$ to commute with the $G$-action.
\end{Lemma}

\begin{proof}
Let $k$ be the order of the cyclic group $\Gamma$.
The bundle $\Nu$, being an equivariant complex vector bundle,
decomposes into a direct sum of equivariant complex vector bundles,
$\Nu = \oplus_{l=1} ^ {k-1} \Nu_l$,
where $z \in \Gamma$, thought of as a $k$th root of unity,
acts on the fibers of $\Nu_l$ by complex multiplication by $z^l$. 
(There is no $\Nu_0$ because
$\Gamma$ has no fixed points outside of the zero section.) 
We define the $S^1$-action on $\Nu_l$ by letting $z \in S^1$ 
act by complex multiplication by $z^l$.
\end{proof}

\begin{proof}[Proof of Lemma \ref{step}]
Denote by $\Gamma$ the $S^1$-stabilizer of $X$.
Let the quotient circle $S^1/\Gamma$ act effectively 
on a closed disk $D^2$ by rotations.
This quotient acts on the stratum $X$ freely, and $X$
is cobordant to zero via the associated disk bundle:
\begin{equation} \labell{X X'}
 X = \partial (X \times_{S^1/\Gamma} D^2).
\end{equation}

Let $\Nu$ be the normal bundle to $X$ in $Z$.
Lemma \ref{extend} gives a fiberwise $S^1$-action on $\Nu$.
The action of $S^1/\Gamma$ on $X$ lifts to an action on $\Nu$
via the anti-diagonal embedding,
$S^1 / \Gamma \mapsto S^1 \times_\Gamma S^1$, given by
$a \mapsto [a,a^{-1}]$, followed by the action of $S^1 \times_\Gamma S^1$
in which the first $S^1$ acts as before and the second $S^1$ acts by the
fiberwise action.
The vector bundle $\Nu \times_{S^1/\Gamma} D^2$ 
over $X \times_{S^1/\Gamma} D^2$ extends 
the normal bundle $\Nu$ over $X$ to the cobordism \eqref{X X'}.

Let $U$ be a tubular neighborhood of $X$ in $Z$ whose closure,
$\ol{U}$, is identified with the closed unit disk bundle in $\Nu$
(with respect to some invariant fiberwise metric).  Then
\begin{equation} \labell{with corners}
  \ol{U} \times_{S^1/\Gamma} D^2,
\end{equation}
is a manifold with boundary and corners.
One piece of its boundary is $\ol{U} \times_{S^1/\Gamma} \partial D^2$,
which we identify with $\ol{U}$. The other boundary piece is
$S(\Nu) \times_{S^1/\Gamma} D^2$. The two pieces intersect along a corner,
$S(\Nu) \times_{S^1/\Gamma} \partial D^2$, which can be identified 
with $S(\Nu)$. The manifold $Z'$ is obtained by the surgery which 
replaces the first boundary piece, $\ol{U}$, by the second. 
To see that this results in a cobordant manifold, thicken $Z$
into $Z \times [0,1]$, and attach \eqref{with corners} to
$\ol{U} \times \{1\}$.
After smoothing the corner at $\partial \ol{U} \times \{1\}$, 
we get a cobordism between $Z$ and $Z'$.
(Also see an alternative proof below.)

The $S^1$-stabilizer of a point of the form $[p,0]$
in $\ol{U} \times_{S^1/\Gamma} D^2$ is equal to the 
stabilizer of $p$ with respect to the fiberwise $S^1$-action.
(Indeed, for $a \in S^1$, let us denote the original $S^1$ action 
on $p \in \Nu$ by $a \cdot p$ and the fiberwise action by $p \star a$.
Then $a$ belongs to the stabilizer of $[p,z]$ if and only if
there exists $b \in S^1$ such that 
$(a \cdot p , z) = (b \cdot p \star b\inv , b\inv \cdot z)$.
Setting $z=0$, it is easy to check that $[a \cdot p , 0] = [p \star a , 0]$,
so if $p \star a = p$, then $a$ is in the stabilizer of $[p,0]$.
Conversely, suppose $(p \star a , 0) = ( b \cdot p \star b\inv , 0)$.
Then $b \cdot p$ and $p$ must be in the same fiber of $\Nu$.
Therefore $b \in \Gamma$, so $b \cdot p \star b\inv = p$,
and so $p \star a = p$.)
This stabilizer is smaller than $\Gamma$ for points in
$X' = \partial \ol{U} \times_{S^1/\Gamma} \{0\}$.
It is all of $S^1$ when $p$ is in $X$.

One can easily verify that
an additional $G$-action on $Z$ extends to the cobordism and hence
to $Z'$, and so does an equivariant stable complex structure.
Similarly, an orientation on $Z$ induces an orientation on the cobordism
and on $Z'$.

Given an abstract moment map $\Psi$ on $Z$, we extend it 
to an abstract moment map $\tPsi$ on $\ol{U} \times_{S^1/\Gamma} D^2$ 
by setting 
$$ \tPsi ([a,z]) = \Psi (g(z)a)$$
where $g \colon D^2 \to \R$ is $S^1$-invariant, equal to $1$ for $z$
near the boundary $\partial D^2$, and equal to $0$ for $z$ near the
origin.
\end{proof}

Although the method used in the above proof is standard in topology,
we provide an alternative construction, in which we need not smooth corners.

\begin{proof}[Alternative proof of Lemma \ref{step}]
Recall that $U$ denotes an open neighborhood of $X$ in $Z$,
whose closure, $\ol{U}$, is identified with the unit disk bundle 
of the normal bundle of $X$ in $Z$. 
The boundary $\partial \ol{U}$
is identified with the unit sphere bundle, $S(\Nu)$.
Recall that $D^2$ denotes the closed unit disk in $\R^2$.  
We can identify 
\begin{equation} \labell{radial}
 \ol{U} \ssminus X \cong S(\Nu) \times \R_+ 
 \quad \quad \text{and} \quad \quad
 D^2 \ssminus \{ 0 \} \cong \partial D^2 \times \R_+,
\end{equation}
where $\R_+ = [0,\infty)$, by identifying $(a,t)$ on the right
with $a/(1+t)$ on the left in both cases;
similarly, setting $\R_{>0} = (0,\infty)$,
\begin{equation} \labell{radial2}
   U \ssminus X \cong S(\Nu) \times \R_{>0}.
\end{equation}
Let $\tN$ be the manifold $\ol{U} \times_{S^1/\Gamma} D^2$
minus its corners. Consider the subset
$\tX = (X \times_{S^1/\Gamma} D^2) 
        \cup (\ol{U} \times_{S^1/\Gamma} \{ 0 \} )$.
Its complement is 
$$\begin{array}{ccl}
\tN \ssminus \tX  & = &  
( (\ol{U} \ssminus X) \times_{S^1/\Gamma} (D^2 \ssminus \{ 0 \}))
\ssminus (\text{\emph{its corners}}), \\
 & = & (S(\Nu) \times \R_+) \times_{S^1/\Gamma} (\partial D^2 \times \R_+)
 \ssminus (\text{\emph{its corners}})
 \quad \text{by \eqref{radial}} \\
 & = &  (S(\Nu) \times_{S^1/\Gamma} \partial D^2) 
	\times (\R_+^2 \ssminus \{ 0 \} ) \\ 
 & = &  S(\Nu) \times (\R_{>0} \times [0,\pi/2])  
\quad \text{in polar coordinates}.
\end{array}$$
Combining with \eqref{radial}, we get
\begin{equation} \labell{common}
 \tN \ssminus \tX = S(\Nu) \times \R_{>0} \times [0,\pi/2] 
 = (U \ssminus X) \times [0,\pi/2].
\end{equation}
We get the cobordism $\tZ$ by gluing $\tN$ with the trivial cobordism
$ (Z \ssminus X) \times [0,\pi/2]$ along their common open set
\eqref{common}.
As before, one can routinely check that the requirements
of Lemma \ref{step} are satisfied.
\end{proof}

\section{Proper cobordisms result in no new relations between compact
cobordism classes}
\labell{sec:compact-versus-proper}

In Section \ref{sec:cutting} we showed that if two compact $G$-manifolds
with proper moment maps are properly cobordant, the manifolds are also 
compactly cobordant, but through an orbifold. Using the surgery of
Section \ref{sec:surgery}, we can obtain a genuine non-singular 
compact cobordism. The following result implies Theorem \ref{thm:intro}:

\begin{Theorem} \labell{thm:stable-complex}
If two compact $G$-equivariant (oriented) stable complex manifolds 
with abstract moment maps are properly cobordant, they are also 
compactly cobordant.
\end{Theorem}

\begin{proof}
To prove the theorem we repeat the cutting procedure from the proof
of Theorem \ref{orbifold}, but desingularize the resulting orbifold
each time after cutting, as described in Section \ref{sec:surgery}.
To be more precise, the Gusein-Zade surgery is applied before taking
the $S^1$-quotient. The resulting manifold inherits an equivariant
stable complex structure, which allows one to repeat this process.
\end{proof}

\begin{Remark}
Recall that a stable complex manifold is orientable, not oriented
(see Appendix \ref{sec:stable-complex}).
Theorem \ref{thm:stable-complex} holds with or without an orientation
on $M$ having been fixed. 
\end{Remark}

\section{Null-cobordant manifolds}
\labell{sec:null}

We now state an easy but important corollary of Lemma \ref{desingularize}.  
In the following theorem, cobordism is
understood as that of stable complex oriented $S^1$-manifolds, with
proper abstract moment maps in the non-compact case.

\begin{Theorem} \labell{thm:GZ}
Let $N$ be a stable complex (oriented) $S^1$-manifold with a locally
free $S^1$-action and a proper abstract moment map. 
\begin{enumerate}
\item 
Then $N$ is properly cobordant to a manifold with a free $S^1$-action.
\item Assume that $N$ is compact. Then $N$ is compactly cobordant to zero.
\end{enumerate}
The above assertions remain correct when $N$ and the cobordisms
are equipped with an action of a compact group $G$, which commutes with 
the $S^1$-action and preserves the abstract moment map. 
\end{Theorem}

\begin{Remark}
The cobordism of the first assertion has fixed points if the action on $M$ 
is not free. The cobordism referred to in the second assertion always has 
fixed points.
\end{Remark}

\begin{proof}[Proof of Theorem \ref{thm:GZ}]
Lemma \ref{desingularize} provides a cobordism
between $N$ and a manifold $N'$ with a free circle action.
This is further cobordant to the empty set
via the associated disk bundle $N' \times_{S^1} D^2$.
\end{proof}

In the rest of this section we show that one can replace in 
Theorem \ref{thm:GZ} the assumption that $N$ is equivariantly stable
complex with certain restrictions on the stabilizer groups which occur
in $N$.

The stable complex structure was used in the surgery to ensure the following

\begin{Condition} \labell{bundle complex}
The normal bundle to any orbit type stratum
admits an equivariant (fiberwise) complex structure.
\end{Condition}

We can assume Condition \ref{bundle complex} to begin with,
instead of working with stable complex manifolds.
To deduce Lemma \ref{desingularize}, we must show that Condition 
\ref{bundle complex} persists through the surgery of Lemma \ref{step}:

\begin{Lemma} \labell{persists}
Let $Z$ be a manifold with an $S^1$-action
that has no fixed points and that satisfies Condition \ref{bundle complex}.
Let $X \subseteq Z$ be a minimal orbit type stratum.
Let $Z'$ and $\tZ$ be the manifold and the cobordism obtained from $Z$
by performing along $X$ the surgery of Lemma \ref{step}.
Then $Z'$ and $\tZ$ also satisfy Condition \ref{bundle complex}. 
\end{Lemma}

\begin{proof}
Every stratum $Y$ in $Z$ other than $X$ naturally extends 
to a stratum $\tilde{Y}$ in $\tZ$ by taking
$Y \times [0,1]$ outside $\ol{U} \times_{S^1/\Gamma} D^2$
and attaching $(Y \cap \ol{U}) \times_{S^1/\Gamma} (D^2 \ssminus \{ 0 \})$.
The invariant complex structure on the normal bundle of $Y$
induces an invariant complex structure on the normal bundle of $\tilde{Y}$.

A stratum in $\tZ$ which does not come from a stratum in $Z$
must have the form  $Y \times_{S^1/\Gamma} \{ 0 \}$
where $Y \subseteq \ol{U}$ is an orbit type stratum for the fiberwise action.
The fibration $\ol{U} \to X$ restricts to a fibration $Y \to X$.
The normal bundle to $Y$ in $\ol{U}$ is the pullback under this fibration
of a complex sub-bundle of the complex vector bundle $\Nu$ over $X$.
This together with the complex structure on $\R^2 = \C$ gives an invariant 
complex structure on the normal bundle to $Y \times_{S^1/\Gamma} \{ 0 \}$
in $\ol{U} \times_{S^1/\Gamma} D^2$.
\end{proof}

This implies the following variant on Lemma \ref{desingularize}:

\begin{Lemma} \labell{desingularize-A}
Let $Z$ be an $S^1$-manifold that has no fixed points
and that satisfies Condition \ref{bundle complex}.
Let $Y \subset Z$ be the set of points where the action is not free.
Then one can obtain a manifold with a free $S^1$-action
by performing an equivariant surgery on $Z$ along $Y$.
The resulting manifold, $Z'$, is equivariantly cobordant to $Z$.
\end{Lemma}

\begin{proof}
Repeat the proof of Lemma \ref{desingularize},
that is, perform iterations of the surgery of Lemma \ref{step}
to get rid of all the non-trivial stabilizers.
Lemma \ref{step} can be applied because of Condition \ref{bundle complex}
and Lemma \ref{persists}.
\end{proof}

The following equivalent form of Condition \ref{bundle complex} 
will be useful:

\begin{Condition} \labell{restriction}
If an orbit type stratum $X$ is contained in the closure of
another orbit type stratum $Y$, 
and if $[\Stab(X):\Stab(Y)]  = 2$, 
then the normal bundle of $X$ in the closure of $Y$ 
admits an equivariant fiberwise complex structure.
\end{Condition}

\begin{Lemma} \labell{conditions}
A manifold with an $S^1$-action satisfies Condition \ref{bundle complex}
if it satisfies Condition \ref{restriction}.
\end{Lemma}

\begin{proof}
Denote by $\Gamma$ the stabilizer of $X$.  The fiberwise action of $\Gamma$
on $\Nu$ uniquely decomposes into a direct sum, where each summand is a
multiple of a non-trivial irreducible real representation of $\Gamma$.
This defines a decomposition of $\Nu$ into a direct sum of real vector bundles.
It suffices to find an invariant complex structure on each summand.
On a summand on which all the elements of $\Gamma$ act by multiplication
by $\pm 1$, the existence of an equivariant complex structure is guaranteed
by Condition \ref{restriction}.

Fix any other summand, $\Nu'$. 
Let $u$ be a generator of $\Gamma$. There exists an $\alpha\in\C$ such that on
every fiber of $\Nu'$ the representation of $\Gamma$ is real equivalent
to the representation of $\Gamma$ on $\C^n$ with $u$ acting
by scalar multiplication by $\alpha$. Fix such an $\alpha$.
Since, by assumption, not all the elements of $\Gamma$ act by multiplication 
by $\pm 1$, we must have that $\alpha\neq\pm 1$. Then
there exist (unique) real numbers $a$ 
and $b$ such that $\sqrt{-1}=a\alpha+b\alpha^2$. The fiberwise operator
$a u+b u^2$ defines an invariant complex structure on $\Nu'$.

Suppose $\Gamma = \Z_k$ with $k$ odd.
The complex number $\alpha$ of the previous paragraph must satisfy
$\alpha^k=1$, hence, $\alpha = u^l$ for some $0 < l < k$. 
If $l$ is odd, the fiberwise circle action on this summand, $\Nu'$, has
stabilizers $\{ e \}$ and $\Z_l$, which are of odd order.
If $l$ is even, by flipping the complex structure on $\Nu'$
we get that $u$ acts as multiplication by $u^{-l} = u^{k-l}$,
with $k-l$ an odd number strictly between $0$ and $k$.
Finally, the stabilizer groups for the fiberwise $S^1$-action
are subgroups of the stabilizer groups for the restriction of this
action to the summands, $\Nu'$, hence are all of odd order.
\end{proof}

Condition \ref{restriction} automatically holds if there are no two
orbit type strata with stabilizers $\Z_m$ and $\Z_{2m}$.
We therefore have the following variant on Theorem \ref{thm:GZ}: 

\begin{Theorem} \labell{thm:GZ-A}
Let $N$ be an (oriented) $S^1$-manifold with a locally
free $S^1$-action and a proper abstract moment map. 
Suppose that the set of subgroups of $S^1$ which occur as stabilizers 
of points in $N$ does not contain both $\Z_m$ and $\Z_{2m}$ 
for any $m$. 
\begin{enumerate}
\item 
Then $N$ is properly cobordant to a manifold with a free $S^1$-action.
\item 
Assume that $N$ is compact. Then $N$ is compactly cobordant to zero.
\end{enumerate}
The above assertions remain correct when $N$ and the cobordisms
are equipped with an action of a compact group $G$, which commutes with 
the $S^1$-action and preserves the abstract moment map. 
\end{Theorem}

\begin{proof}
The condition on the stabilizers in $N$ automatically implies Condition 
\ref{restriction}. By Lemma \ref{conditions}, this implies Condition
\ref{bundle complex}. Lemma \ref{desingularize-A} provides an equivariant
cobordism between $N$ and a manifold $N'$ with a free $S^1$-action. 
This is further cobordant to the empty set
via the associated disk bundle $N' \times_{S^1} D^2$.
\end{proof}

\section{Restrictions on stabilizers}

Conner and Floyd, as well as Gusein-Zade, considered (compact)
cobordisms in which all the stabilizer subgroups are constrained to lie
in a pre-assigned collection of subgroups. (See, e.g., \cite{CF:book},
\cite{CF:complex}, and \cite{GZ:short}.)  

Let $\calA$ be a collection of subgroups of $G$.  A
\emph{$(G,\calA)$-manifold} is by definition a $G$-manifold such that all
stabilizers of the $G$-action belong to $\calA$.  Define the
\emph{$(G,\calA)$-compact cobordism group} to be the group (under
disjoint unions) of all compact $(G,\calA)$-manifolds with abstract
moment maps modulo compact $(G,\calA)$-cobordisms with abstract moment
maps, and define the \emph{$(G,\calA)$-proper cobordism group} to be
the group of all (possibly non-compact) $(G,\calA)$-manifolds with
proper abstract moment maps modulo $(G,\calA)$-cobordisms with proper
abstract moment maps.

For instance, the cobordisms of locally free actions
are obtained by setting $\calA$ to be the collection of
all discrete subgroups of $G$. The cobordisms of quasi-free
actions result from taking the collection of connected subgroups of $G$
as $\calA$.

We expect that an analogue of Theorem \ref{thm:stable-complex}
would be true when instead of working with equivariant stable complex 
manifolds we work with $(G,\calA)$-manifolds, for certain collections
$\calA$. Namely, we expect that under appropriate conditions on $\calA$,
the $(G,\calA)$ compact cobordism group would inject into the $(G,\calA)$
proper cobordism group. For instance, this is true for quasi-free circle
actions, when $G = S^1$ and $\calA = \{ S^1 , \{ 1 \} \}$, by Remark
\ref{quasi free}. More generally, when $G=S^1$, it is reasonable to demand 
that $\calA$ does not contain both $\Z_m$ and $\Z_{2m}$ for any $m$. 
This condition automatically implies Condition \ref{restriction},
hence (by Lemma \ref{conditions}) Condition \ref{bundle complex},
hence it allows the surgery of Lemma \ref{desingularize}. However, 
this condition in itself is not sufficient; the combination of cutting
and surgery might introduce new stabilizers:

\begin{Example}
Suppose that $\calA$ does not contain the group $G=S^1$.
Then the $(G,\calA)$ compact cobordism group does not inject into the
$(G,\calA)$ proper cobordism group. For instance, $M=S^1$ 
is properly cobordant to the empty set through the proper
cobordism $S^1 \times \R_+$, on which the action is free.
However, $S^1$ is not locally-freely compactly cobordant to the empty
set; the quotient of such a cobordism would be a compact one-dimensional
manifold whose boundary consists of a single point, and this cannot happen.
\end{Example}

\begin{Remark}
Compact equivariant cobordisms have been studied in the past in the presence 
of equivariant stable complex structures or with certain restrictions 
on the stabilizer subgroups.

For instance, Conner and Floyd, in their study of cobordisms 
equivariant with respect to cyclic groups,
treated involutions separately from maps of odd period in \cite{CF:book},
and worked with an equivariant stable complex structure in \cite{CF:complex}.

Gusein-Zade, in his study of $S^1$-equivariant oriented cobordisms
in \cite{GZ:short} and \cite{GZ:long}, assumed either 
the existence of a stable complex structure,
or that all finite stabilizers are of odd order.
In fact, his reason for these assumptions is the same as ours -- 
to allow the surgery of section \ref{sec:surgery}.
He also mentioned the case where the possibility for $\calA$ to consist
of $S^1$ together with all finite cyclic groups $\Z_m$
for which the decomposition of $m$ into prime factors contains an even 
number of $2$'s. This restriction, too, implies Condition \ref{restriction}.
\end{Remark}

\appendix
\section{Stable complex structures}
\labell{sec:stable-complex}
A (tangential) \emph{stable complex structure} on a manifold $M$ 
is represented 
by a complex structure on the fibers of the real vector bundle 
$TM \oplus \R^k$ for some $k$;
two complex structures, on the fibers of $TM \oplus \R^k$ and of
$TM \oplus \R^l$, give rise to the same stable complex structure
if and only if there exist $a$ and $b$ such that the complex vector bundles
$TM \oplus \R^k \oplus \C^a$ and $TM \oplus \R^l \oplus \C^b$
are isomorphic. To make sense of stable complex cobordisms,
notice that a stable complex structure on $M$ induces one on $\partial M$.

In algebraic topology, a stable complex bundle is usually defined
on the stable normal bundle to $M$ (see, e.g., \cite{stong} and
\cite{rudyak}). In the non-equivariant setting
this is equivalent to defining a tangential stable complex structure.
In the equivariant case, the situation is different and there is no
one-to-one correspondence between these structures. (See \cite{may},
p.~337 for details.) For example, many a result of Section 
\ref{sec:compact-versus-proper} would not be true for manifolds with
stable complex structures on the normal bundle. For this reason,
in the present paper as, e.g., in \cite{may},
we define stable complex structures on $G$-manifolds to be on 
stable tangent bundles. Sometimes, (tangential)
stable complex structures are also referred to as stable almost 
complex structures \cite{may,rudyak} or weakly almost complex structures
\cite{BH}.

Stable complex structures arise, for example, from complex
structures, almost complex structures, or symplectic structures.
However, in general, a stable complex structure does not have to be
associated with an almost complex structure.

A stable complex manifold is orientable, but not canonically oriented.
For instance, on $\C$, $i$ and $-i$ induce the same stable complex
structure, but opposite orientations. The reason is that $T\C$ 
with either of these structures is equivalent to the trivial complex line
bundle over $\C$.
To see that a
stable complex manifold is orientable, it suffices to notice that the
sum of $TM$ and a trivial bundle is a complex (hence, oriented) bundle.
By fixing an orientation on the trivial bundle, we get an orientation on
$TM$.
Because of this, in the theory of stable complex cobordisms (also referred
to as unitary or complex cobordisms) an orientation 
is often fixed in addition to a stable complex structure.
(See, e.g., \cite{stong}, \cite{rudyak}, and \cite{may}.) 

An \emph{equivariant stable complex structure} on a $G$-manifold 
is an equivalence class of invariant complex structures on 
$TM \oplus \R^k$. The $G$-action on this bundle is standard on the first
term and trivial on the second. The equivalence relation is
$G$-equivariant complex vector bundle equivalence.

\begin{Remark}
Because a stable complex structure is, by definition, an equivalence
class, a group action on a manifold does not induce a group action 
on the set of stable complex structures. Therefore, one cannot say that
a given stable complex structure, in our sense, is invariant with respect
to the group action. This is the reason why an equivariant stable complex
structure is defined as an equivalence class of invariant complex structures
on $TM \oplus \R^k$.
\end{Remark}

Note that two equivariantly
distinct equivariant stable complex structures can be equal
as non-equivariant stable complex structures.
For example, take $M=\C$ and let $G=S^1$ act by rotations.  Then $i$ and
$-i$ define the same stable complex structure but distinct equivariant
stable complex structures.

In this paper we use the following facts:
\begin{enumerate}
\item
A stable complex structure on a manifold naturally induces a stable
complex structure on any regular level set of a smooth function on this
manifold.
\item
An equivariant stable complex structure on a manifold with a free action
of a compact group naturally induces a stable complex structure on the
quotient manifold.
\item
When $M$ has an equivariant stable complex structure,
each orbit type stratum inherits a stable complex structure, 
and its normal bundle is naturally a complex vector bundle
on which the group acts by complex bundle automorphism. (This would
not be true for ``normal'' equivariant stable complex structures discussed
above.) 
\end{enumerate}

There are different notions that are related to the notion of
a stable complex structure used in this paper. For example,
weakly complex structures from \cite{CF:complex,CF:relation} 
are similar to our stable complex
structures, but the equivalence relation is finer; two representatives,
$J_1$ and $J_2$, on $TM \oplus \R^k$, give rise to the same weakly complex
structure
if and only if $J_1$ and $J_2$ are homotopic after stabilization.
Every weakly complex structure 
gives rise to a stable complex structure and to
an orientation on the manifold. This shows that the notion of a 
weakly complex structure 
is not equivalent to our notion of a stable complex structure.

\subsection*{Acknowledgement}
We would like to thank Piotr Mikrut for his valuable
remarks and comments on the proof of Lemma \ref{desingularize}.

\end{document}